\documentclass{amsart}
\usepackage{amsmath}
\usepackage{amsthm}
\usepackage{amssymb}
\usepackage{hyperref}

\newtheorem{thm}{Theorem}

\newtheorem{prop}{Proposition}

\newcommand{\be}{\begin{equation}}
 \newcommand{\ee}{\end{equation}}
 \newcommand{\ba}{\begin{array}{l}}
 \newcommand{\ea}{\end{array}}
 \newcommand{\pa}{\partial}
\newcommand{\na}{\nabla}

 \newcommand{\la}{\label}
\newcommand{\fr}{\frac}

\newcommand{\dm}{{\mbox{div}}_m}

\newcommand{\nam}{\nabla_m}
\newcommand{\dam}{\mbox{div}_m}

\newcommand{\nax}{\nabla_x}

\newcommand{\dx}{{\mbox{div}}_x}

\newcommand{\Rr}{{\mathbb{R}}}
\newcommand{\e}{\epsilon}

\newcommand{\Det}{\mbox{Det}}
\newcommand{\Tr}{\mbox{Tr}}

\begin{document}

\title
{Note on Global Regularity for 2D Oldroyd-B Fluids with Diffusive Stress}

\author{Peter Constantin}

\address{Department of Mathematics, Princeton University\\Fine Hall, Washington Road\\Princeton, New Jersey 08540}

\email{const@math.princeton.edu}

\author{Markus Kliegl}
\address{Department of Mathematics, The University of Chicago\\5734 S. University Ave\\Chicago, Illinois 60637}

\email{kliegl@princeton.edu}

\begin{abstract} We prove global regularity of solutions of Oldroyd-B equations
in 2 spatial dimensions with spatial diffusion of the polymeric stresses.
\end{abstract}


\maketitle

\section{Introduction}
We discuss a simple model of polymeric fluid, a complex fluid that is comprised of a solvent, which is an incompressible Newtonian fluid, and a dilute suspension of polymeric matter in it. The complex fluid occupies a region in physical space $\Rr^d$, and in this paper $d=2$. In models (\cite{doied}, \cite{otbook}),  complicated objects are described by retaining very few degrees of freedom, and the polymers are represented by end-to-end vectors  $m\in M=\Rr^d$. In general, the configuration space $M$ could be much more complicated. We consider as starting point the kinetic descriptions of the particles. The state 
of the particles is determined by a measure $f(m)dm$, i.e. by a measure which is absolutely continuous with respect to the measure $dm$.  The equilibrium measure is obtained by minimizing a modified free energy
\be
{\mathcal{E}}[f] = \int_M\left(\log f + U_1(m) +\fr{1}{2}U_2[f](m)\right)fdm
\la{free}
\ee
In this description $f\log f$ represents a thermal effect, $U_1(m)$ is a resident potential that describes constraints and  
\be
U_2[f](m) = \int_M k(m,p) f(p)dp
\la{u2}
\ee
is a screening (or excluded volume) potential. If the kernel $k(m,p)$ is symmetric and Lipschitz then the resulting extremal equation (the Onsager equation),
\be
f= Z^{-1}e^{-U[f]}\la{ons}
\ee
has always solutions, at least if $M$ is compact (\cite{cons},\cite{consz}).
 Here
\be
U[f] = U_1 + U_2[f]\la{u}
\ee
 and $Z>0$ is a normalizing constant. The modified free energy is not convex, and the Onsager equation has rigorously proven phase transitions (when temperature is lowered) and multiple solutions coexist. The limit of high intensity, (or low temperature) is determined by a popularity contest (\cite{consz}): a selection mechanism operates by which states which can be easily altered are rejected in favor of states whose character is very much like that of their neighbors.
The relaxation mechanisms by which the minimum solutions of Onsager equations are approached is given by the kinetic equation
\be
\pa_t f = \e{\mbox{div}}_m\left (f\na_m\left(\fr{\delta {\mathcal{E}}}{\delta f}\right)\right)
\la{kin}
\ee
with $\e$ a positive constant quantifying inter-particle diffusivity. The equation can be written as
\be
\pa_t f = \e\Delta_m f + \e{\mbox{div}}_m\left(f\na_m(U[f])\right)\la{fu}
\ee
and is a nonlinear Fokker-Planck equation. The presence of the nonlinearity shows that the problem of deriving this equation from some underlying stochastic system is not trivial.
The equation (\ref{kin}) has $\mathcal E$ as a Lyapunov functional. The time derivative of this Lyapunov functional is non-positive, vanishes at solutions of the Onsager equation, and if $M$ is connected, only there:
\be
\fr{d}{dt}{\mathcal{E}} = - {\mathcal {D}}
\la{ed}
\ee
where
\be
{\mathcal {D}} = \e\int_M f\left|\nam\left(\log f + U[f]\right)\right |^2dm\la{mathcald{d}}
\ee
The presence of fluid introduces new degrees of freedom, due to the ambient space, a subset of $\Rr^d$ or $\mathbb T^d$. The particle density acquires $x$-dependence: $f= f(x,m,t)$. Now the kinetic evolution depends on the fluid's velocity $u(x,t)$. This velocity is a function of the macroscopic variables $x,t$ alone. The equation ceases to be a gradient equation. The fluid velocity introduces a drift both in space and in particle configuration space. The simplest case in which such an effect is immediately observed is the case in which a shear transports an elastic rod. The nonlinear Fokker-Planck equation is then
\be
\pa_t f + u\cdot\nax f +(\nax u)m\cdot\nam f = \e\Delta_m f + \e\dam(f\nam U[f])
\la{fum}
\ee
Note that if $u$ is divergence-free in $x$ then $(\nax u)m$ is divergence-free in $m\in \Rr^d$. The ensemble of particles might suffer collectively the effect of being in a strained environment. This can be modeled by allowing the resident potential $U_1$ to depend on $x,t$ as well as on $m$. The resident potential does not depend on the state $f$. The self-interaction potential $U_2[f]$ depends on the macroscopic variables only because $f$ does, but the kernel $k$ is derived from purely microscopic information. For more general configuration spaces we have an added drift
\be
W(x,m,t) = \nax u(x,t) : c(m)\la{w}
\ee
that depends linearly on the spatial gradient of velocity, and gives a vector field on $M$. In components:
$$
W^{\alpha}(x,m,t) = (\pa_j u^i(x,t))c^{j, \alpha}_i(m)
$$
The coefficients $c(m)$ are smooth.
The nonlinear Fokker-Planck equation is
\be
\left(\pa_t + u\cdot\nax\right )f + \dm(W f) = \e\Delta_m f + \e\dm(f\nam U[f])
\la{nlfp}
\ee
In this equation the coefficients $c$ represent a ``macro-micro'' interaction that introduces a particle drift from the macroscopic drift. The potential kernel $k$ embodies a ``micro-micro'' interaction, while the configuration space $M$ and the resident potential $U_1$ represent geometric and kinematic constraints. The effect that particles might have on the solvent is a more mysterious matter. At this level of description, in order to be self-consistent, this effect can only be obtained by averaging out the microscopic variables. The effect is embodied in an added stress matrix $\sigma(x,t)$ and it is a ``micro-macro'' interaction, the macroscopic effect of microscopic insertions. 
An energetic principle (\cite{c}) states that this ``micro-macro'' interaction, is such that the coupled system is dissipative.  It turns out that this principle is sufficient to provide formulas for the ``micro-macro'' interaction even in non-dilute cases.
In the known examples, this principle leads to familiar rules of determining the added polymeric stress from the micro-micro and the macro-micro interactions (\cite{lebris}). In (\cite{crem}) this principle is described for non-dilute situations.
In general, the fluid velocity solves the incompressible (unit density) Navier-Stokes equation
\be
\pa_t u + u\cdot\nax u -\nu\Delta_x u +\nax p = K\dx\sigma
\la{nse}
\ee
with $K$ a positive constant with units of velocity squared ($cm^2/sec^2$),
\be
\nax\cdot u =0
\la{divu}
\ee
and $\nu>0$ the kinematic viscosity. The dimensionless matrix $\sigma = (\sigma_{ij}(x,t))$ represents the added stress tensor. We note that only the non-isotropic part of $\sigma$ enters the equation, because changing $\sigma$ to
$\sigma + \lambda {\mathbb I}$ results in the same equation, with a modified pressure. The coupled system (\ref{nlfp}, \ref{nse}) is required to dissipate the sum of the kinetic energy
and free energy:
\be
\fr{d}{dt}\int_{\Rr^d} \left[\fr{1}{2}|u(x,t)|^2 + K{\mathcal {E}}[f](x,t)\right]dx + \int_{\Rr^d} \left[\nu|\nax u(x,t)|^2 + K{\mathcal {D}}(x,t)\right]^2dx \le 0
\la{requ}
\ee
In (\cite{c}, \cite{crem}) it is shown that this requirement  (for all initial data) is satisfied if, and only if
\be
\sigma(x,t) =\int_M \left [ c\cdot\nam U - \dm c\right]fdm 
\la{sig}
\ee  
and
\be
(\pa_t + u\cdot\nax) U_1\le 0.
\la{dtu}
\ee
Recall that $c$ is a matrix of vector fields so $c\cdot\nam U$ and $\dm c$  are linear operators in  $\Rr^d$.
In the case of the Oldroyd B equation, $U_2=0$ and 
\be
U_1(m) = \frac{|m|^2}{2R^2}.
\la{um}
\ee 
$R$ is positive, $M=\Rr^d$, and
\be
W(x,t) =\nax u(x,t)m
\la{wo}
\ee  
i.e. $c(m)^{j,\alpha}_i = \delta^{\alpha}_im^{j}$. The formula (\ref{sig}) gives
\be
\sigma(x,t) = \int_{\Rr^d}(m\otimes \na_m U)fdm
\la{sigconc}
\ee
(Note that $\dm c =\pa_{\alpha}(\delta^{\alpha}_im^j) =0$.)

In the Oldroyd B case $\sigma$ obeys a closed equation, and this is easily obtained by multiplying (\ref{fum}) by $(m\otimes m)/R^2$ and integrating. The result is:
\be
(\pa_t + u\cdot\nax)\sigma = (\nax u)\sigma + \sigma (\nax u)^T - 2k\sigma + 2k\rho {\mathbb I} \la{sigreq}
\ee
where the damping frequency is
\be
k = \fr{\e}{R^2}
\ee
and
\be
\rho(x,t) = \int_{\Rr^d}f(x,m,t)dm
\la{rhoo}
\ee
Note that from (\ref{fum}) it follows that
\be
(\pa_t + u\cdot\nax)\rho =0.
\la{rhoeqq}
\ee

The mathematical literature on complex fluids is growing, and we cannot give here a complete account. Early work (\cite{saut}, \cite{renardy}) established local existence results for Oldroyd-B and FENE type equations. The Oldroyd-B equations are exact closures of linear Fokker-Planck equations. The resident potential is a harmonic potential that is independent of macroscopic variables and the
self-interaction kernel $k$ vanishes. The ``micro-macro'' interaction gives $\sigma$ as a second moment of $f$ and then a self-contained equation for $\sigma$ follows from the linear Fokker-Planck equation. This is the only known instance when the Fokker-Planck equation yields a self-contained equation for $\sigma$. The FENE equation is a model in which the resident potential is infinity at a finite extension value, prohibiting the particles from extending beyond it. 
Global existence of weak solutions in the presence of spatial diffusion of the polymers was proved in a sequence of papers, the most recent of which is (\cite{suli}). Global existence of weak solutions via propagation of compactness was proved under the corotational assumption (\cite{lionsmasmoudi}, \cite{lionsmasmoudi1}) and, very recently, for the full FENE model (\cite{mas1}). There is no such result for the Oldroyd B model. The global existence of smooth solutions for small data for Oldroyd B-type models was established in (\cite{leiz}, \cite{fanghuasmall}) and for FENE in (\cite{masfene}). 
Global existence of smooth solutions for large data in 2D was established for Smoluchowski equations on compact manifolds (\cite{c-smo}, \cite{cftz}, \cite{cm}, \cite{cs}, \cite{cs1}, \cite{otto}). Global regularity for large data in the FENE case, under the corotational assumption was proved in (\cite{fanghua},\cite{mas}). An approach based on Lagrangian particle dynamics was developed in (\cite{fanghuachun}). Sufficient conditions for regularity in terms of bounds on the added stress tensor were established in (\cite{chemm}, \cite{kmt}) and further refined in (\cite{leim}). Numerical evidence for singularities was provided in (\cite{thomasses}). The paper (\cite{cws}) proved global existence for small data with large gradients for Oldroyd-B. A regularization obtained by allowing the ``spring constant'' in the harmonic potential to depend on the local rate of strain of the fluid was obtained also in (\cite{cws}). Regularity for diffusive Oldroyd-B equations in 2D for large data were obtained in the creeping flow regime in (\cite{crem}). More general models will be presented in (\cite{ckli}).

\section{Diffusive Oldroyd B: a priori bounds}
We have only limited success in the proof of regularity for large data.
The case of (\ref{nlfp}) coupled with the 2D Navier-Stokes equations, when $M$ is a compact Riemannian manifold was addressed in (\cite{cm}) using Fourier analysis techniques and in (\cite{cs}, \cite{cs1}) using physical space techniques. If the configuration space of the particles is not compact, then the only results of global regularity for all data are for modified equations: the corotational case when  2D Navier-Stokes equations are coupled via the antisymmetric part of the gradient to linear FENE equations (\cite{masfene}), the case of Oldroyd B coupled with 2D Navier-Stokes when the resident potential responds to excessive rate of strain in the fluid (\cite{cws}), and the case in which the polymeric stress is allowed to diffuse in space (\cite{crem}, in the creeping flow regime). 
The proofs of global existence for the compact configuration space $M$ are based on the fact that $\sigma\in L^{\infty}(dx dt)$, which easily follows from definitions. The $2D$ Navier Stokes equations forced by the divergence of a bounded stress have unique weak solutions that are very well behaved. In particular, they are H\"{o}lder continuous after an initial transient time. This fact was proved in (\cite{cs}). The proof (\cite{cs1}) of global regularity of the coupled system uses then a local existence theorem to address the initial transient time, and to obtain a global-in-time uniform bound for the velocity in a H\"{o}lder space. Then a bootstrap for higher regularity is used to finish the proof. The proof of the analogous result in (\cite{cm}) uses ideas from (\cite{chemm}) to derive estimates in which the regularity is deteriorating in time, but in a manner that is controlled locally uniformly. Assuming a finite time singularity, we reach a contradiction by starting the computation close to the putative singularity and showing that the controlled loss of singularity forces the solution to remain smooth beyond the blow-up time.  The case of $M=\Rr^d$ with spatial diffusion for $\sigma$, for Oldroyd B in a creeping flow regime (coupling to Stokes) in 2D was presented in (\cite{crem}), and the coupling to the full Navier-Stokes equation is in this paper. For the FENE equations (see \cite{suli1} for numerical results and set-up) (\cite{masfene}) proves local well-posedness in general, and global well posedness for the corotational case.
The corotational case is one in which the full gradient $\nabla u$ is replaced by its anti-symmetric part in (\ref{w}). In that case, if the resident potential is radial, then the correct $\sigma$ given by the relation (\ref{sig}) vanishes. If that is the case there is no added stress to the fluid, and coupling with (\ref{sigconc}) is not energetically balanced.

The 2D Oldroyd B system is the equation (\ref{sigreq}) coupled with (\ref{nse}), (\ref{divu}), (\ref{rhoeqq}).  We consider the variables
\be
\left\{
\ba
a(x,t) = \frac{1}{2}\left ( \sigma^{11}(x,t) - \sigma^{22}(x,t)\right),\\
b(x,t) = \sigma^{12}(x,t) = \sigma^{21}(x,t),\\
c(x,t) = \sigma^{11}(x,t) + \sigma^{22}(x,t) = \Tr\,(\sigma(x,t))
\ea
\la{abc}
\right.
\ee

The matrix $\sigma $ is symmetric and positive by construction, in view of  (\ref{sigconc}), and is given in terms of $a,b,c$ by

\be
\sigma = \left (
\begin{array}{cc}
\fr{c}{2} + a & b\\
b & \fr{c}{2} -a
\end{array}
\right )
\la{sigmabc}
\ee
The positivity of the matrix is equivalent to the positivity of $c$ (which follows from (\ref{sigconc})) and the positivity of the determinant, i.e. 
\be
\fr{c^2}{4}-a^2-b^2>0.
\la{detpos}
\ee
We denote by
\be
\lambda(x,t) = \fr{1}{2}\left(\pa_1 u^1(x,t) -\pa_2u^2(x,t)\right),
\la{lam}
\ee
\be
\mu(x,t) =\fr{1}{2}\left(\pa_1 u^2(x,t) + \pa_2 u^1(x,t)\right)\la{mu}
\ee
and
\be
\omega(x,t) = \pa_1u^2(x,t)-\pa_2u^1(x,t)\la{ome}
\ee
The functions $\lambda$ and $\mu$ represent the rate of strain, $\omega$ the vorticity.
The equations (\ref{sigreq}) can be written as the system
\be
\left\{
\ba
D_t a = -\omega b + \lambda c - 2ka \\
D_t b  = \omega a + \mu c -2k b,\\
D_t c = 4\lambda a + 4\mu b -2kc + 4k\rho\\
D_t\rho = 0
\ea
\right.
\la{abceq}
\ee
We used the notation $D_t = \pa_t + u\cdot\nax$. 
Let us multiply the $c$ equation by $\fr{c}{2}$, the $a$ equation by $2a$, the $b$ equation by $2b$ and subtract the last two from the first. We obtain
\be
D_t\left( \fr{c^2}{4} - a^2-b^2\right) = -4k\left (\fr{c^2}{4}-a^2-b^2\right) + 2k\rho c
\la{deteq}
\ee
This cancellation of nonlinearity is not surprising because
\be
\fr{c^2}{4} - a^2-b^2 = \Det\, (\sigma)
\la{detdet}
\ee
and the determinant is conserved along particle trajectories if $k=0$.
It is well-known, and easy to see that the regularity of the system is decided by whether or not we can bound $c$ in $L^{\infty}$.

We prove global existence in the case of physical space diffusion of $\sigma$.
The system we consider is thus
\be
(\pa_t + u\cdot\nax)\sigma = \sigma(\nax u) + (\nax u)^T\sigma -2k(\sigma-\rho{\mathbb I}) + \kappa\Delta_x\sigma
\la{sigdiff}
\ee
with $\kappa>0$ a spatial diffusivity of the stress, coupled with (\ref{nse}), (\ref{divu}), (\ref{rhoeqq}). In terms of
$a,b,c, \rho$, this is the system
\be
\left \{
\ba
(\pa_t + u\cdot\nax) a = -\omega b + c\lambda -2k a +\kappa\Delta_x a\\
(\pa_t + u\cdot\nax)b = \omega a  +c\mu -2kb +\kappa\Delta_x b\\
(\pa_t + u\cdot\nax) c =4(\lambda a + \mu b) - 2kc + \kappa\Delta_x c + 4k\rho \\
(\pa_t + u\cdot\nax) \rho = 0\\
\ea
\la{abcd}
\right .
\ee
coupled with (\ref{nse}), (\ref{divu}).
We assume we are in two spatial dimensions, that the initial data for $a,b,c, \rho$ are in $L^1(\Rr^2)\cap W^{1,2}(\Rr^2)$. We assume that the initial data for the velocity are divergence-free and belong to $W^{2,2}(\Rr^2)$. We start with a calculation that shows that $c$ is bounded below in time, as long as $(a,b,c)$ solve  an equation of the form (\ref{abcd}) with bounded $\lambda$ and $\mu$. It does not matter how the functions $u$, $\lambda$, $\mu$ and $\omega$ are obtained, nor even if they are related among themselves.
We consider the equation obeyed by
\be
\gamma(x,t) = c-2\sqrt{a^2+b^2}
\la{gamma}
\ee
A brief calculation shows
\be
D_t\gamma = -2\left( k + \fr{\lambda a +\mu b}{\sqrt{a^2+b^2}}\right)\gamma
+ 4k\rho +\kappa\Delta_x c -\fr{2\kappa}{\sqrt{a^2 +b^2}}(a\Delta_x a+ b\Delta_xb)
\la{dtga}
\ee 
Now we have
\be
\Delta_x\left(\sqrt{a^2+b^2}\right) -\fr{a\Delta_xa+ b\Delta_xb}{\sqrt{a^2+b^2}} =
(a^2+b^2)^{-\fr{3}{2}}\sum_{j}(a\pa_j b-b\pa_j a)^2\ge 0
\la{wa}
\ee
so we deduce
\be
D_t\gamma \ge -2\left( k + \fr{\lambda a +\mu b}{\sqrt{a^2+b^2}}\right)\gamma
+\kappa\Delta_x \gamma +  4k\rho.  
\la{dtgam}
\ee
In view of the maximum principle, and the fact that $\rho\ge 0$, it follows that $\gamma(x,t)\ge 0$ if the initial data obey $\gamma(x,0)\ge 0$. If the initial data $\sigma_0$ is nonnegative, i.e. it has nonnegative determinant, $d_0=\Det\sigma_0\ge 0$ and nonegative trace $c_0\ge 0$, then $\gamma_0 = 2\sqrt{a_0^2+b_0^2 +d_0}-2\sqrt{a_0^2+b_0^2}\ge 0$, and then it follows that
\be
c\ge2\sqrt{a^2+b^2}\la{clow}
\ee
for all time. This implies that $c\ge 0$ and that $c^2\ge 4(a^2+b^2)$, so we proved
\begin{prop} \la{p1} A solution of (\ref{sigdiff}) in 2D, advected by a $L^1(0,T; W^{1,\infty}(\Rr^2))$ incompressible velocity, stays nonnegative if it starts nonnegative.
\end{prop}

\noindent{\bf Remark.} This is not an isolated result. In fact any equation of the type
\be
D_t \sigma = M\sigma + \sigma M^T + \kappa \Delta \sigma
\la{genm}
\ee
with $M$ traceless and bounded, preserves positivity. The proof above is special to two dimensions, but the result is true in any number of dimensions, as is easily proved by time-splitting. Indeed, if $\kappa=0$, the equation preserves the determinant, and hence positivity. The heat equation, on the other hand, preserves positivity because its solution is obtained by convolution with a positive function.

We are ready now to obtain good global apriori bounds. 
We consider thus the coupled system
\be
\left\{
\ba
(\pa_t + u\cdot\nax )\sigma = (\nax u)\sigma + \sigma (\nax u)^T -2k(\sigma -\rho {\mathbb I}) +\kappa \Delta_x\sigma\\
(\pa_t + u\cdot\nax )\rho =0\\
(\pa_t + u\cdot\nax -\nu\Delta_x) u +\nax p = K\dx\sigma\\
\nax\cdot u= 0
\ea
\right .
\la{usig}
\ee
in ${\mathbb R^2}$. We recall that $u$ has dimensions of velocity $cm/sec$, $\rho$ and $\sigma$ are dimensionless, $k$ has dimensions of inverse time $sec^{-1}$, $\nu$ and $\kappa$ have dimensions of diffusivity $cm^2/sec$, and $K$ of $cm^2/sec^2$. The equations are dimensionally correct: the first two and the last one have dimensions of inverse time and  the third has dimensions of $cm/sec^2$. The bounds below can be checked for units; this is a way of double-checking that the bounds make sense, and is why we keep the constant physical coefficients  $k, \kappa, \nu, K$, without making nondimensional combinations out of them.

We start by using the energy balance, taking advantage of the fact that $\sigma$ is a positive matrix in view of Proposition {\ref{p1}}.
We take the trace of the first equation, integrate in space and add to the third equation multiplied by $2u$ and integrated. We obtain, after the cancellation of the cubic terms
\be
\fr{d}{dt}\int_{\Rr^2}\left [|u|^2 + K\Tr(\sigma )\right]dx + \int_{\Rr^2}\left [2\nu |\nax u|^2 + 2kK \Tr(\sigma)\right]dx \le 4kK\int_{\Rr^2}\rho_0(x)dx
\la{en}
\ee
where $\rho_0$ is the initial density of particles. We used the fact that $\rho$ is transported by a divergence-free vector field, and so all its spatial integrals are conserved. We obtain thus 
\be
\sup_{0\le t\le T}\left[ \|u\|_{L^2}^2 + K\|\sigma\|_{L^1}\right] + 2\nu\int_0^T\|\nax u\|_{L^2}^2dt \le \|u_0\|_{L^2}^2 +K\|\sigma_0\|_{L^1} +4kKT\|\rho_0\|_{L^1} = R_0
\la{eniq}
\ee
We use the fact that $c\ge 2\sqrt{a^2+b^2}$ and that $\int_{\Rr^2} cdx$ is 
equivalent to the $L^1$ norm of $\sigma$. The right hand side is bounded uniformly, independently of $\nu, \kappa$. We take now the first equation, multiply 
by $\sigma$ and integrate.  We obtain
\be
\fr{d}{2dt}\|\sigma\|_{L^2}^2 + \kappa\|\nax \sigma\|_{L^2}^2 + 2k\|\sigma\|^2_{L^2}\le 2\|\sigma\|_{L^4}^2\|\nax u\|_{L^2} + 2k\|\rho_0\|_{L^2}\|\sigma\|_{L^2}
\la{sigl2}
\ee 
Using the Ladyzhenskaya inequality
\be
\|\sigma\|_{L^4}^2\le C \|\sigma\|_{L^2}\|\nax \sigma\|_{L^2}
\la{lady}
\ee
valid in 2d, we deduce 
\be
\fr{d}{dt}\|\sigma\|_{L^2}^2 +\kappa \|\nax\sigma\|_{L^2}^2 \le
C\left[\kappa^{-1}\|\nax u\|_{L^2}^2\right ]\|\sigma\|_{L^2}^2 + k\|\rho_0\|^2_{L^2}
\la{sigl2i}
\ee
and, by Gronwall:
\be
\sup_{0\le t\le T}\|\sigma\|^2_{L^2} + \kappa\int_0^T\|\nax \sigma\|^2_{L^2}\le
Ce^{\fr{\|u_0\|_{L^2}^2 + K\|\sigma_0\|_{L^1} + 4kKT\|\rho_0\|_{L^1}}{\nu\kappa}}\left[\|\sigma_0\|_{L^2}^2 + kT\|\rho_0\|_{L^2}^2\right] = R_1
\la{exp}
\ee
It follows that $\sigma\in L^{\infty}(dt, L^2(\Rr^2))$ and
$\sigma\in L^2(dt, W^{1,2}(\Rr^2))$ are bounded apriori in terms of the initial data. This used the previous bound (\ref{eniq}), and that used the positivity of $\sigma$. Taking the curl of the third equation in (\ref{usig}) we obtain
\be
\pa_t \omega + u\cdot\nax\omega -\nu\Delta\omega = K\nax^{\perp}\cdot\dx\sigma
\la{omegak}
\ee 
Multiplying by $\omega$ and integrating we obtain, after integration by parts and Young's inequality:
\be
\fr{d}{dt}\|\omega\|^2_{L^2} +\nu\|\nax\omega\|^2_{L^2}\le \fr{K^2}{\nu}\|\nax \sigma\|_{L^2}^2
\la{naxomega}
\ee
Using (\ref{exp}) we have an apriori bound for the time integral of the right-hand side, and so we obtain
\be
\ba                                                                            
\sup_{0\le t\le T}\|\omega\|^2_{L^2} + \nu\int_0^T\|\nax \omega\|^2_{L^2}\\ \le   \fr{CK^2}{\kappa\nu}e^{\fr{\|u_0\|_{L^2}^2 + K\|\sigma_0\|_{L^1} + 4kKT\|\rho_0\|_{L^1}}{\nu\kappa}}\left[\|\sigma_0\|_{L^2}^2 + kT\|\rho_0\|_{L^2}^2\right] + \|\omega_0\|^2_{L^2} \\ 
= R_2
\ea
\la{omegab}
\ee
Now we take the first equation in (\ref{usig}) multiply by $-\Delta_x\sigma$ and integrate.
We obtain
\be
\fr{d}{2dt}\|\nax \sigma\|_{L^2}^2 + \kappa\|\Delta_x\sigma\|^2_{L^2} \le
\|\nax u\|_{L^{2}}\|\nax \sigma\|_{L^4}^2 + 2\left [\|\sigma\|_{L^4}\|\nax u\|_{L^4} + k\|\rho\|_{L^2}\right]\|\Delta_x \sigma\|_{L^2} 
\la{sigma1}
\ee
The first term in the right-hand side comes from the advective term after integration by parts and use of incompressibility. We bound it using Ladyzhenskaya's inequality and Young's inequality. We arrive at
\be
\fr{d}{dt}\|\nax \sigma\|_{L^2}^2 + \kappa \|\Delta_x \sigma\|_{L^2}^2 \le
C{\fr{\|\nax u\|_{L^2}^2}{\kappa}}\|\nax\sigma\|_{L^2}^2 + \fr{C}{\kappa}\|\sigma\|_{L^4}^2\|\nax u\|_{L^4}^2 + \fr{Ck^2}{\kappa}\|\rho_0\|_{L^2}^2
\la{naxsig1}
\ee
From Ladyzhenskaya's inequality, (\ref{lady}) applied to $\sigma $ and to $\nax u$, from (\ref{exp}) and (\ref{omegab}) and usual elliptic estimates relating $\nax u$ in $L^p$ to $\omega$ in the same space, we see that
\be
\fr{1}{\kappa}\int_0^T\|\sigma\|_{L^4}^2\|\nax u\|_{L^4}^2\le B
\la{B}
\ee
where the bound $B$ is an explicit expression in terms of the initial data, coefficients and $T$:
\be
B = \fr{C}{\kappa\sqrt{\kappa\nu}}R_1R_2
\la{Bou}
\ee
where $R_1$, $R_2$ are the right-hand sides of (\ref{exp}), (\ref{omegab}). Note that $B$ is dimensionless. It follows then that
\be
\ba
\sup_{0\le t\le T}\|\nax\sigma\|_{L^2}^2 + \kappa\int_0^T\|\Delta\sigma\|_{L^2}^2dt \\\le Ce^{\fr{\|u_0\|_{L^2}^2 + K\|\sigma_0\|_{L^1} + 4kKT\|\rho_0\|_{L^1}}{\nu\kappa}}\left[\|\nax\sigma_0\|_{L^2}^2 + B + \fr{k^2T}{\kappa}\|\rho_0\|_{L^2}^2\right]\\
=R_3
\ea
\la{naxsigb}
\ee
Armed with this, we return to (\ref{omegak}), multiply by $\Delta_x\omega$ and integrate in space. We obtain, after usual manipulations
\be
\fr{d}{dt}\|\nax \omega\|^2_{L^2} + \nu\|\Delta \omega\|^2_{L^2} \le C\fr{\|u\|_{L^{\infty}}^2}{\nu}\|\nax \omega\|_{L^2}^2 + C\fr{K^2}{\nu}\|\Delta_x\sigma\|_{L^2}^2
\la{naxom1}
\ee
We use the interpolation inequality
\be
\|u\|_{L^{\infty}}^2 \le C \|u\|_{L^2}\|\Delta u\|_{L^2}
\la{interfty}
\ee
valid in $\Rr^2$ and thus, using (\ref{eniq}), (\ref{omegab}) and (\ref{naxsigb}) we obtain
\be
\ba
\sup_{0\le t\le T}\|\nax \omega\|_{L^2}^2 + \nu \int_0^T\|\Delta\omega\|_{L^2}^2dt \\
\le Ce^{C\nu^{-\fr{3}{2}}\sqrt{TR_0R_2}}\left[\|\nax\omega_0\|^2_{L^2} + \fr{K^2}{\nu\kappa}R_3\right] = R_4
\ea
\la{deltom}
\ee
We differentiate (\ref{rhoeqq})
\be
\pa_t (\nax\rho) + u\cdot\nax(\nax\rho) = - (\nax u)^T\nax\rho
\la{naxrhoeq}
\ee
We use the interpolation inequality
\be
\|\nax u\|_{L^{\infty}}^2\le C \|\nax u\|_{L^2}\|\Delta \nax u\|_{L^2}
\la{naxuint}
\ee
to deduce that
\be
\sup_{0\le t\le T}\|\rho\|_{W^{1,2}(\Rr^2)} \le e^{\nu^{-\fr{1}{4}}R_2^{\fr{1}{4}}T^{\fr{3}{4}}R_4^{\fr{1}{4}}}\|\rho_0\|_{W^{1,2}(\Rr^2)}= R_5
\la{naxrhob}
\ee

We proved:
\begin{prop}{\la{p2}} Let $u,p, \sigma,\rho$ be a strong solution of (\ref{usig}) on $[0,T]$, with $\sigma_0 \in L^1(\Rr^2)\cap W^{1,2}(\Rr^2)$ a nonnegative symmetric matrix, $\rho_0\in L^1(\Rr^2)\cap W^{1,2}(\Rr^2)$ nonnegative and $u_0\in W^{2,2}(\Rr^2)$ divergence-free. 

Then $\rho\in L^{\infty}(0, T, L^1(\Rr^2)\cap W^{1,2}(\Rr^2))$, 
$u\in L^{\infty}(0,T, W^{2,2}(\Rr^2))\cap L^2(0,T, W^{3,2}(\Rr^2))$ and
$\sigma\in L^{\infty}(0,T, L^1(\Rr^2)\cap W^{1,2}(\Rr^2))\cap L^2(0,T, W^{2,2}(\Rr^2))$ obey the bounds (\ref{eniq}), (\ref{exp}), (\ref{omegab}),  (\ref{naxsigb}), (\ref{deltom}), and (\ref{naxrhob}).
\end{prop}
The bounds above used the equations and are fragile: they use the positivity of the matrix $\sigma$. There are several methods of approximating solutions, but these (Galerkin approximation, mild formulation and Picard iteration) do not preserve positivity. The simplest proof of global regularity is therefore done then by continuation. We show (by any method) that a
strong solution exists and is unique for a short time. Because it solves the nonlinear equation, if $\sigma$ starts positive, it stays positive by Proposition \ref{p1}, and then, by the a priori bounds, the solution has controlled size at the end of the short time existence. Existence and uniqueness of strong solutions then allows to conclude that the maximal time of existence is infinite.

\section{Diffusive Oldroyd B: Local and Global Regular Solutions}
\begin{prop}{\la{p3}} Let $u_0\in W^{2,2}(\Rr^2)$ be a divergence-free vector. Let
$\sigma_0\in W^{1,2}(\Rr^2)$ be a symmetric matrix, and let $\rho_0\in L^1(\Rr^2)\cap W^{1,2}(\Rr^2)$. There exists a time $T_0>0$ depending on the norms of
$u_0, \sigma_0, \rho_0$ and a unique strong solution 
\be
\ba
u\in L^{\infty}(0,T_0, {\mathbb P}W^{2,2}(\Rr^2))\cap L^2(0,T_0, {\mathbb P}W^{3,2}(\Rr^2))\\
\sigma \in  L^{\infty}(0,T_0, W^{1,2}(\Rr^2))\cap L^2(0,T_0, W^{2,2}(\Rr^2))\\
\rho\in L^{\infty}(0, T_0, L^1(\Rr^2)\cap W^{1,2}(\Rr^2))
\ea
\la{belongs}
\ee
of (\ref{usig}) with initial data $u_0,\sigma_0, \rho_0$. 
\end{prop} 
Above ${\mathbb P}$ is the Leray-Hodge projector on divergence-free vectors. 
We sketch a proof via an explicit scheme, a Picard iteration for the mild formulation of the equations. We take the Banach space defined by the right hand sides of (\ref{belongs}) with the corresponding norm $\|(u,\sigma,\rho)\|_B$ (the sum of the respective norms). 
\be 
\|(u,\sigma,\rho)\|_B = \|u\|_X + \|\sigma\|_Y + \|\rho\|_{Z}
\la{norms}
\nonumber
\ee
where $X = L^{\infty}(0, T_0, {\mathbb P}W^{2,2}(\Rr^2))\cap L^2(0,T_0, {\mathbb P}W^{3,2}(\Rr^2))$,
\be
\|u\|_X = \|u\|_{L^{\infty}(0,T_0, W^{2,2}(\Rr^2))} + \|u\|_{L^2(0,T_0, W^{3,2}(\Rr^2))}
\la{normux}
\nonumber
\ee
$Y= L^{\infty}(0, T_0, W^{1,2}(\Rr^2))\cap L^2(0,T, W^{2,2}(\Rr^2))$,
\be
\|\sigma\|_Y = \|\sigma\|_{L^{\infty}(0,T_0, W^{1,2}(\Rr^2))} + \|\sigma\|_{L^2(0,T_0, W^{2,2}(\Rr^2))}
\la{norsigy}
\nonumber
\ee
and
$Z= L^{\infty}(0,T_0, L^1(\Rr^2) \cap W^{1, 2}(\Rr^2))$ with norm

\be
\|\rho\|_Z = \sup_{0\le t\le T_0}\left(\|\rho (t)\|_{L^1(\Rr^2)} + \|\rho(t)\|_{W^{1,2}(\Rr^2)}\right)
\nonumber
\ee
We set up  a fixed point equation $U = F(U)$ in $B$ for the triplet $U= (u,\sigma,\rho)$, where $F(U)= (u^{new}, \sigma^{new}, \rho^{new})$ given by
\be
\ba
u^{new}(t) = e^{\nu t\Delta}u_0 + Q_1(u, u) + L_1(\sigma)\\
\sigma^{new}(t) = e^{(\kappa\Delta -2k)t}\sigma_0 + Q_2(u, \sigma)+ L_2(\rho)\\
\rho^{new} (t) = N(u), \; {\mbox{where}}\; (\pa_t + u\cdot\nax)\rho^{new} = 0, \quad \rho^{new}(x,0) = \rho_0(x)
\ea
\la{set}
\ee
where
\be
Q_1(u,v)(t) = -\int_0^te^{\nu(t-s)\Delta}{\mathbb P}(u(s)\cdot\nax v(s))ds,
\la{q1}
\ee 
\be
L_1(\sigma)(t) = K\int_0^te^{\nu(t-s)\Delta}{\mathbb P}(\dx\sigma(s))ds
\la{l1}
\ee
\be
Q_2(u,\sigma)(t) = \int_0^te^{(t-s)(\kappa\Delta-2k)}(-u(s)\cdot\nax\sigma(s) + (\nax u(s))\sigma(s) + \sigma(s)(\nax u(s))^T)ds
\la{q2}
\ee
and
\be
L_2(\rho)(t) = 2k\int_0^te^{(t-s)(\kappa\Delta-2k)}\rho(s)\mathbb I ds
\la{l2}
\ee
We check that
\be
\ba
\|Q_1(u,v)\|_X\le \epsilon \|u\|_X\|v\|_X\\
\|L_1(\sigma)\|_X\le C\|\sigma\|_Y\\
\|Q_2(u,\sigma)\|_Y \le \epsilon\|u\|_X\|\sigma\|_Y\\
\|L_2(\rho)\|_Y \le \epsilon\|\rho\|_{L^{\infty}(0,T_0, L^{2}(\Rr^2))}\\
\|N(u)\|_{Z}\le C\|\rho_0\|_{L^1(\Rr^2)\cap W^{1,2}(\Rr^2)}e^{C\|u\|_X}\\
\|N(u)-N(v)\|_{L^1\cap L^2}\le \epsilon\|u-v\|_Xe^{C(\|u\|_X + \|v\|_X)}
\ea
\la{ineqs}
\ee
where  $\epsilon$ can be made as small as we wish by choosing $T_0$ small enough (while keeping everything else fixed) and $C$ is a fixed constant. Obviously $Q_1$ and $Q_2$ are bilinear, $L_1$ and $L_2$ are linear, and the map $N$ is nonlinear.
The inequalities (\ref{ineqs}) require little effort to check, and they involve straightforward inequalities including the Ladyzhenskaya inequality and maximal regularity of the heat equation. For instance, the function
\be 
q_1(x,t) = Q_1(u,v)(x,t)\la{qone}
\ee
obeys
\be
\pa_tq_1 -\nu\Delta_xq_1 = {\mathbb P}(u\cdot \nax v), \quad q_1(x,0) = 0
\la{q1eq}
\ee
and the bound is easily obtained by energy methods, via
\be
\fr{d}{2dt}\|\Delta q_1\|_{L^2}^2 + \nu\|\nax\Delta q_1\|_{L^2}^2
\le \|\nax(u\cdot\nax v)\|_{L^2}\|\nax\Delta q_1\|_{L^2},
\la{qoneb}
\nonumber
\ee
the Ladyzhenskaya inequality for the term involving $\nax u\nax v$ and
the interpolation inequality (\ref{interfty}) for the term involving $u \cdot \nax (\nax v)$, to yield:
\be
\ba
\sup_{0\le t\le T}\|\Delta q_1\|_{L^{2}(\Rr^2)}^2 + \nu\int_0^{T_0}\|\nax \Delta q_1\|_{L^{2}(\Rr^2)}^2dt \\
\le \fr{C}{\nu}\int_0^{T_0}[\|\nax u\|_{L^2}\|\Delta u\|_{L^2}\|\nax v\|_{L^2}\|\Delta v\|_{L^2} + \|u\|_{L^2}\|\Delta u\|_{L^2}\|\Delta v\|_{L^2}^2]dt\\
\le CT_0\|u\|_{X}^2\|v\|_{X}^2
\ea
\la{qonebo}
\ee
Next, $l_1(x,t) = L_1(\sigma)(x,t)$ obeys
\be
\pa_t l_1 -\nu\Delta l_1 = K{\mathbb P}(\dx\sigma), \quad l_1(x,0) = 0.
\la{l1eq}
\ee
Applying $\Delta$ and multiplying by $\Delta l_1$ we obtain easily the desired estimate
\be
\sup_{0\le t\le T_0}\|\Delta l_1\|_{L^2}^2 + \nu\int_0^{T_0}\|\nax \Delta l_1\|_{L^2}^2 \le \fr{CK^2}{\nu}\int_0^{T_0}\|\Delta \sigma\|_{L^2}^2 dt
\la{loneb}
\ee
The estimate for $q_2 = Q_2(u, \sigma)$ is straightforward and is easily obtained by multiplying the equation
\be
\pa_t q_2 -\kappa\Delta q_2 +2k q_2= -u\cdot\nax\sigma + (\nax u)\sigma + \sigma(\nax u)^T
\la{q2eq}
\ee
by $\Delta q_2$ and integrating. The estimate for $l_2 = L_2(\rho)$ is obtained in the same fashion. The bound for $N(u)$ is simply obtained by differentiating. The bound for $N(u)-N(v)$ in $L^1\cap L^2$ uses the a priori bound for individual $N(u)$ and $N(v)$.

We set  up the iteration $U^{n+1}= F(U^{n})$. By choosing $\epsilon$ small enough it is then easy to see that the sequence $U^n$ is bounded in $B$ and converges exponentially fast. Let us verify first that we can prove by induction that there exist constants $A, \Gamma, D$ and $\epsilon$ such that
\be
\ba
\|u^n\|_X \le A,\\
\|\sigma ^n\|_Y\le \Gamma,\\
\|\rho ^n\|_Z\le D
\ea
\la{abd}
\ee
By induction we have
\be
\ba
\|u^{n+1}\|_X \le A_0 + \epsilon A^2 + C\Gamma,\\
\|\sigma^{n+1}\|_Y \le \Gamma_0 + \epsilon A\Gamma + \epsilon R_0, \\
\|\rho^{n+1}\|_Z \le R_1e^{CA}
\ea
\nonumber
\ee
where $A_0, \Gamma_0, R_0, R_1$ depend only of the norms of the initial data.
We can choose for instance $A =2(A_0 + C\Gamma)$ and solve first the
induction inequality
\be
\Gamma \ge  \Gamma_0 + 2\epsilon(A_0 + C\Gamma)\Gamma +\epsilon R_0
\nonumber
\ee
by choosing $\epsilon$ small enough. Then, by choosing $\epsilon$ smaller we guarantee that $\epsilon A <\frac{1}{2}$ which implies that
$$
A\ge A_0 + \epsilon A^2 + C\Gamma
$$
Then $D$ is chosen in function of $A$.
The induction succeeds because $\|\rho^n\|_{L^{\infty}(0,T_0, L^2(\Rr^2))}\le R_0$ is bounded a priori. Once we know that the sequence is bounded in $B$, then we can obtain exponential convergence of $u^n$ in $X$, $\sigma^n$ in $Y$ and $\rho^n$
in $L^{1}\cap L^2$.

\begin{thm}{\la{thm1}} Let $u_0\in W^{2,2}(\Rr^2)$ be divergence-free. Let $\rho_0\in L^1(\Rr^2)\cap W^{1,2}(\Rr^2)$ be nonnegative and let $\sigma_0\in L^1(\Rr^2)\cap W^{1,2}(\Rr^2)$ be a nonnegative matrix. Let $T>0$ be arbitrary. There exists a unique strong solution of (\ref{usig}) satisfying the bounds  (\ref{eniq}), (\ref{exp}), (\ref{omegab}),  (\ref{naxsigb}), (\ref{deltom}), and (\ref{naxrhob}).
\end{thm}
The proof follows along the lines described above: by Proposition \ref{p3}
a strong solution exists for a short time interval $[0,T_0]$. We consider the maximal interval of existence: $T_1= \sup T_0\le T$ such that the solution is strong on $[0,T_0]$. By Proposition \ref{p1}
the matrix $\sigma$ remains positive. By Proposition \ref{p2} the solution obeys a priori bounds that are uniform on $[0,T]$. (The a priori bounds are bounded uniformly for $T_0\le T$).  Then, it must be that $T_1=T$, because otherwise, by Proposition \ref{p3} we could extend the solution beyond $T_1$, contradicting its definition.

\begin{thm}\la{thm2} If $u_0\in W^{m+1, 2}(\Rr^2)$ is divergence-free,
$\sigma_0\in L^1(\Rr^2)\cap W^{m,2}(\Rr^2)$ is a nonnegative symmetric matrix and $\rho_0\in L^1(\Rr^2)\cap W^{m,2}(\Rr^2)$ is nonnegative, $m\ge 1$,  then the solution obeys 
\be
\ba
u\in L^{\infty}(0,T, W^{m+1,2}(\Rr^2))\cap L^2(0,T, W^{m+2,2}(\Rr^2)),\\ 
\sigma\in L^{\infty}(0,T, L^1(\Rr^2)\cap W^{m,2}(\Rr^2))\cap L^2(0,T, W^{m+1, 2}(\Rr^2),\\ 
\rho\in L^{\infty}(0,T, L^1(\Rr^2)\cap W^{m,2}(\Rr^2)).
\ea
\nonumber
\ee
\end{thm}
{\bf Remark.} The higher regularity assumption on the initial data is needed because $\rho$ does not obey a diffusive equation and is not constant. If $\rho_0 =1$ then $\rho =1$ for all time. In that case, it follows from parabolic regularization that   $u,\sigma\in C^{\infty}(t_0, T,\Rr^2)$ for positive $t_0$ and the initial data at $t=0$ can be much less regular.
\section{Conclusion}
We proved global regularity for the Oldroyd B system coupled with the full incompressible Navier-Stokes system in 2D when the polymeric stress is diffusing in space. This kind of model was considered for weak solutions in several studies (\cite{suli}) and was used in numerical investigations (\cite{thomasses}, \cite{to}). Our results use the positivity of the stress matrix, which is preserved under diffusive evolution.
\subsection*{Acknowledgment}PC's research was partially supported by NSF-DMS grant 0804380.

\end{document}